\documentstyle[amscd,amssymb,verbatim,epsf,11pt]{amsart}

\def\centereps#1#2#3{\centerline{\epsfxsize#1\epsfysize#2\epsfbox{#3}}}

\begin{document}
\title[Eigenfunctions on a Stadium]{Eigenfunctions on a Stadium Associated with 
Avoided Crossings of Energy Levels}
\author{B. Neuberger, J. W. Neuberger, D. W. Noid}
\email{jwn@@unt.edu}
\email {dwn@@ornl.gov}
\maketitle
\pagenumbering{arabic}
\section {introduction}
   The most primitive form of a stadium is a circle, say of radius $r$. 
If we separate the left and right halves of the circle slightly and connect 
these two semi circles at the top and bottom by straight lines we have 
a classical stadium (racetrack). Let $a$ be the horizontal distance from 
the center of the figure to either the left or right  semi circle; that is, 
the distance of each `straightaway' is $2a$.

    By the 1970's researchers were using numerical studies to examine the
stochastic behavior of quantum mechanical systems corresponding to chaotic
and quasi-periodic behavior in classical mechanics. In an unpublished letter 
written at the Oak Ridge National Laboratories in the early 1980's, \cite{n},
the authors studied a Hamiltonian as a free particle confined in a stadion 
(stadium) boundary. They developed a FORTRAN code to compute the 10 smallest 
eigenvalues of even-even symmetry for the Schroedinger equation for this 
system and plotted (by hand!) line graphs of the eigenvalues versus $a/r$ 
for $r=1$. The correlation plot showed many avoided crossings. Word was that
Wigner proved these graphs could not cross when the corresponding
eigenfunctions had the same symmetry; but after the intervening years a
firm reference evades us. Suggested earlier references include \cite{a} - 
\cite{e}, and \cite{f} for a more recent paper.
Computations in \cite{n} were done only for the upper right quadrant of the 
stadium and the symmetry was forced.

  The present paper is based on the earlier work at Oak Ridge. We have used
a slightly updated version of the old FORTRAN code, this time on the entire
region, again with $r=1$. It is our purpose to indicate certain striking
properties of eigenfunctions near avoided crossings. Plots for these
eigenfunctions were not made at the time \cite{n} was written, nor were
the data retained. 

\newpage
\section{The eigenvalue problem}
   We now consider the eigenvalue problem  
$$ -u_{11}-u_{22}=\lambda u $$
on the interior of the stadium where $u$ is $0$ on the boundary and $\lambda$
represents energy. For each stadium configuration determined by the 
separation distance $a$ the code picks out a designated number of eigenvalues 
$\lambda $ and solves for the values of the associated  eigenfunctions $u$.
%

   Mathematica was used to make surface and contour plots of the eigen
functions. The contour plots provided especially accessible information 
so we concentrated on those. In fact, we visually selected eigenfunctions 
that produced contour graphs that were `even-even', symmetric about both the 
$x$ and $y$ axes. Then we plotted line graphs of the eigenvalues, associated
with selected eigenfunctions, versus $a/r=a$ and examined consecutive pairs
of curves.

   The first figure indicates the third, fourth and fifth eigenvalues 
corresponding to the eigenfunctions of even-even symmetry over a series of values 
$a$ in $[0,2]$. We denote by $\, f_3,  f_4, f_5\, $ the simple curves of those 
eigenvalueswith respect to $a$. The avoided crossings between the pair $(f_3,f_4)$ 
at approximately $.785\, $ and the pair $(f_4,f_5)$, at about $1.51$, are evident 
visually. If one half closes one's eyes it appears that the pairs cross. After a 
closer observation it then appears that the upper curve of each pair continues in 
the same direction its partner was taking before the avoidance, and vice versa.

   The idea of the work here is to examine the eigenfunctions that correlate
with the line curves with particular attention given to contour plots
near the avoided crossings to observe how sensitive the data is to the
separation values $a$. The corresponding pairs of contour plots
($ee3, ee4, ee5$) seem to `swap' characteristics just as the curve pairs 
tended to take over each other's directions.

   Jumping across the narrow gap between $f_4$ and $f_5$ the {\it difference}
between the two eigenfunctions demonstrates a graphical affinity for
the subsequent characteristics of $ee5$ and the previous characteristics
of $ee4$, while the {\it sum} goes the other way (last two contour graphs). 
In this paper we are not drawing any physical conclusions but submit 
some of our graphical results for the reader's inspection.

   Acknowledgement: \hspace {.05cm} We wish to thank A. J. Zettl for his
cheerful help in translating the Wigner-von Neuman paper \cite{a}.
 
\newpage
\vspace*{1.in}
\centereps{12cm}{8cm}{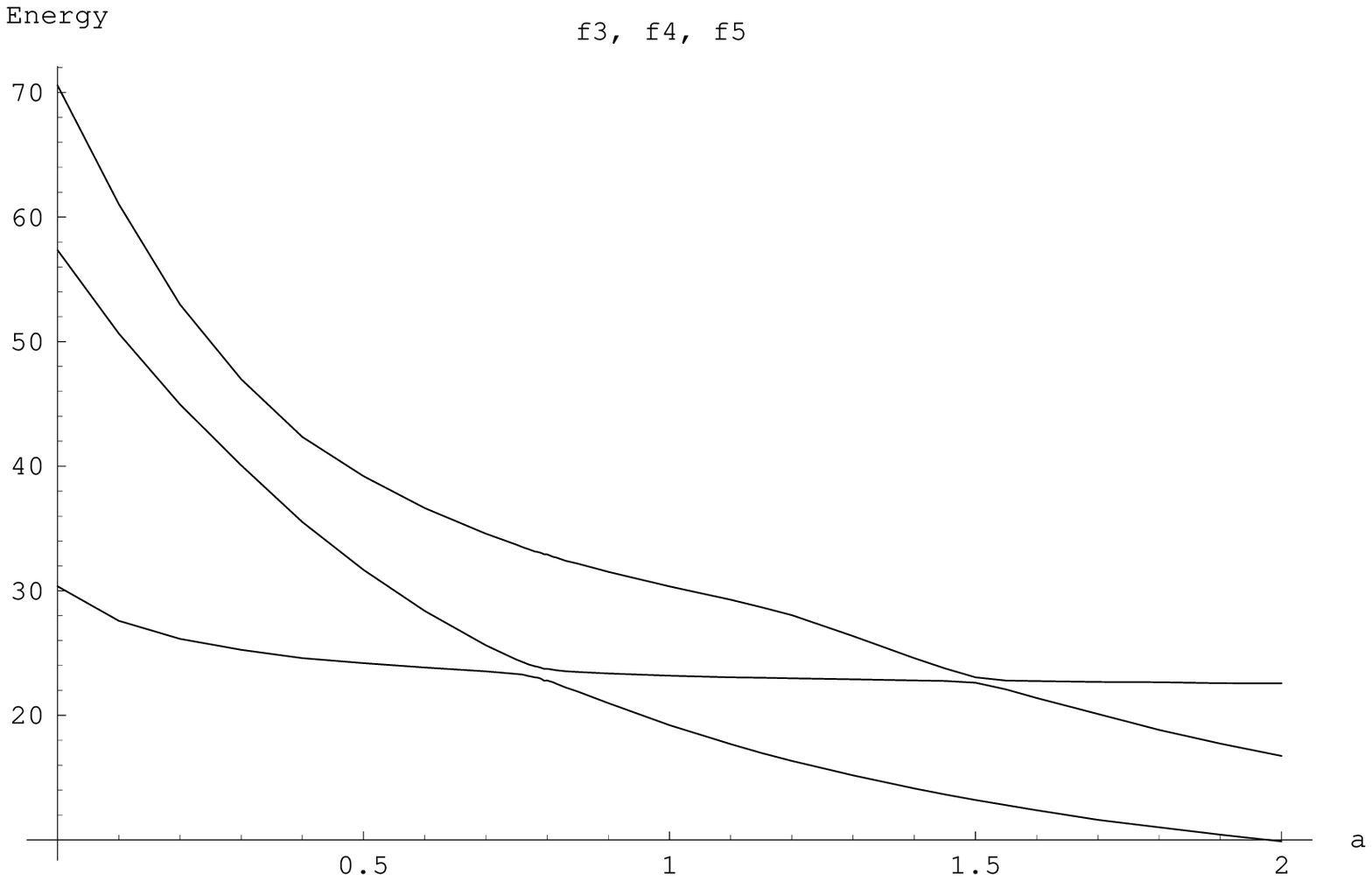}

\newpage{
\vspace*{.1in}
\centereps{9.2cm}{4.cm}{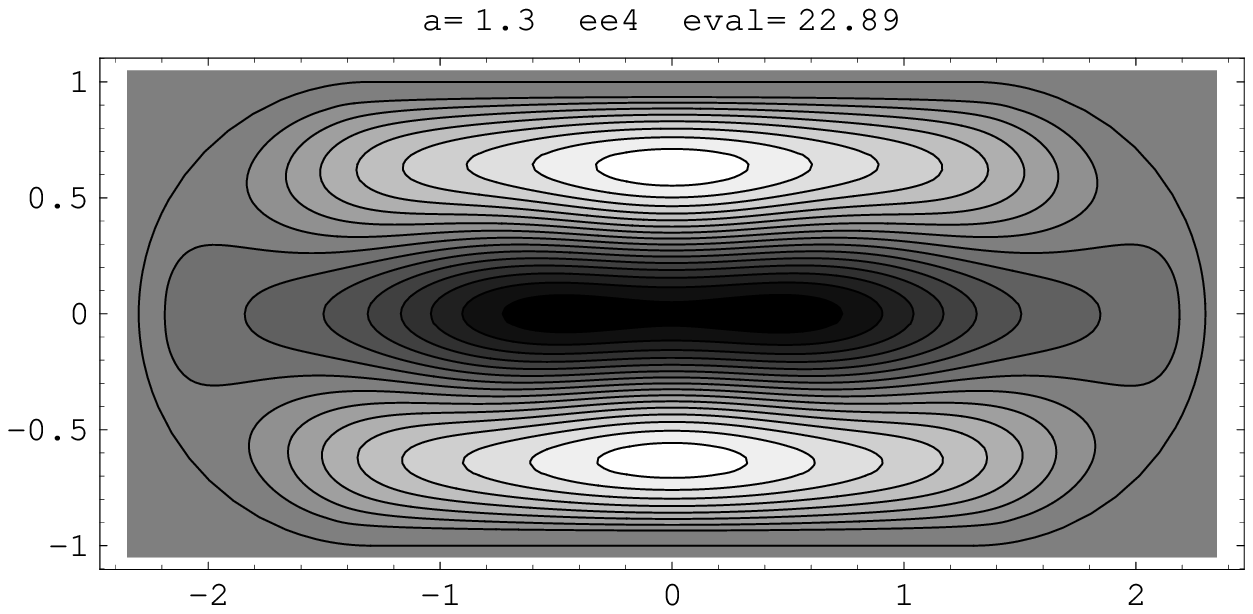}
\centereps{9.2cm}{4.cm}{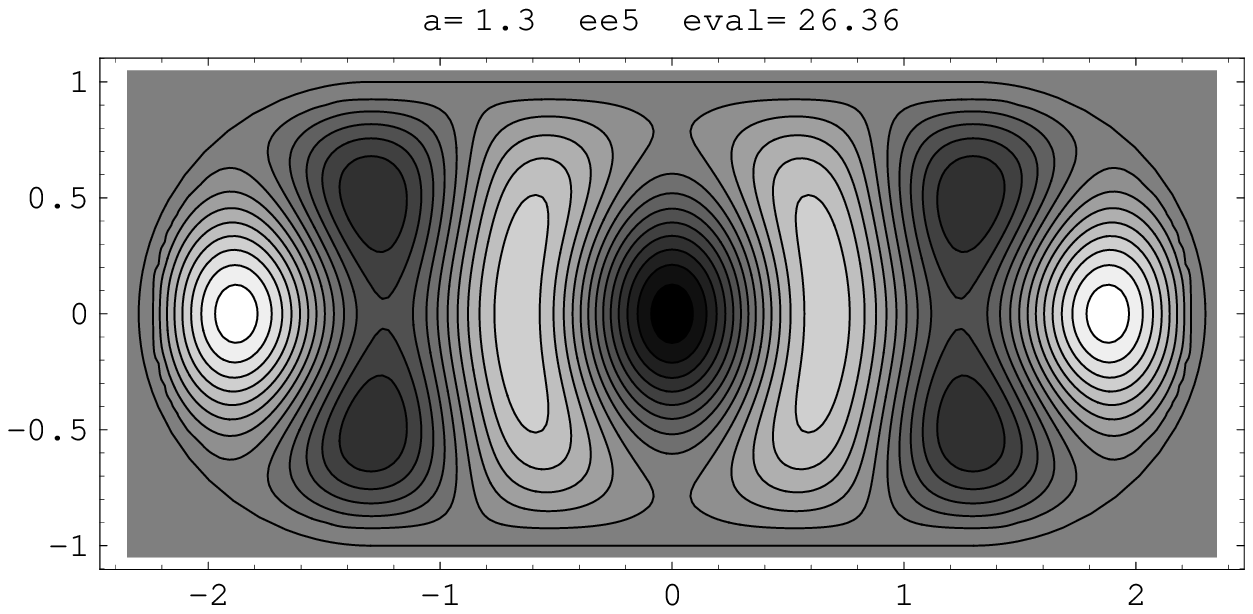}
\centereps{9.8cm}{4.cm}{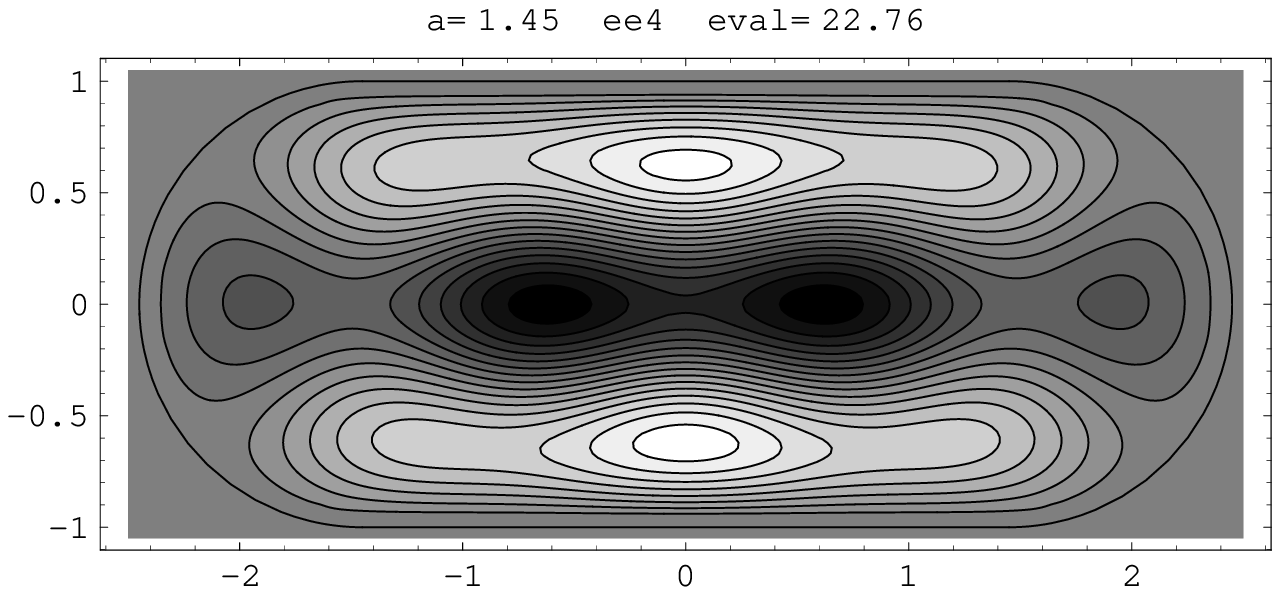}
\centereps{9.8cm}{4.cm}{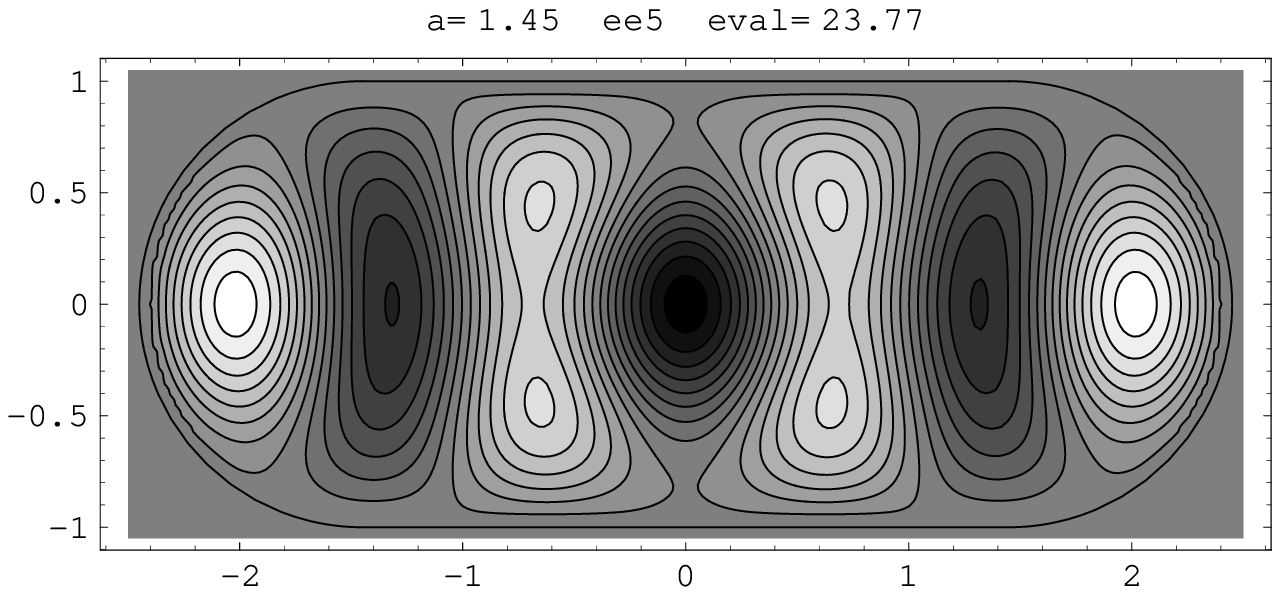}
}

\newpage{
\vspace*{.1in}
\centereps{10.04cm}{4.cm}{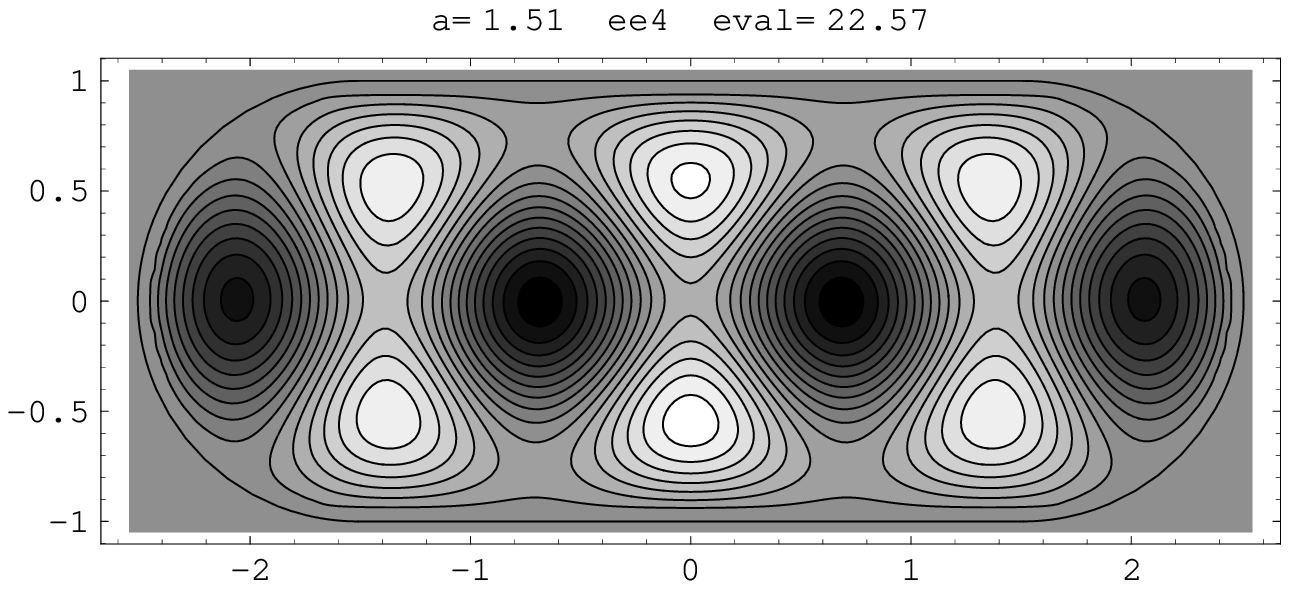}
\centereps{10.04cm}{4.cm}{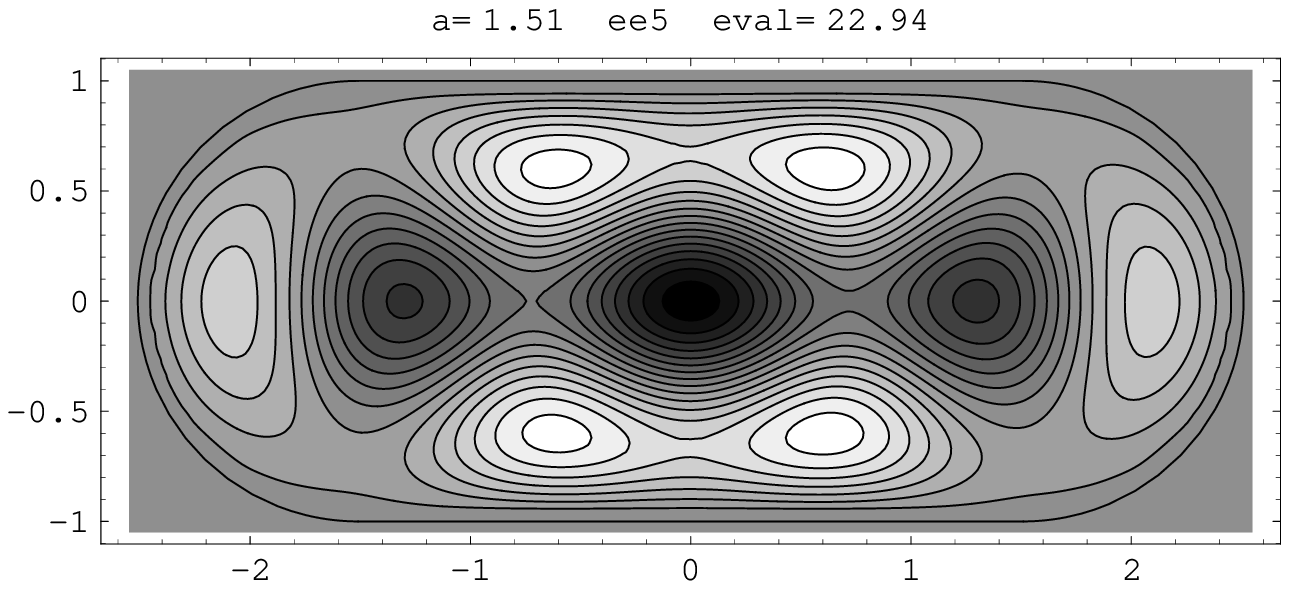}
\centereps{10.2cm}{4.cm}{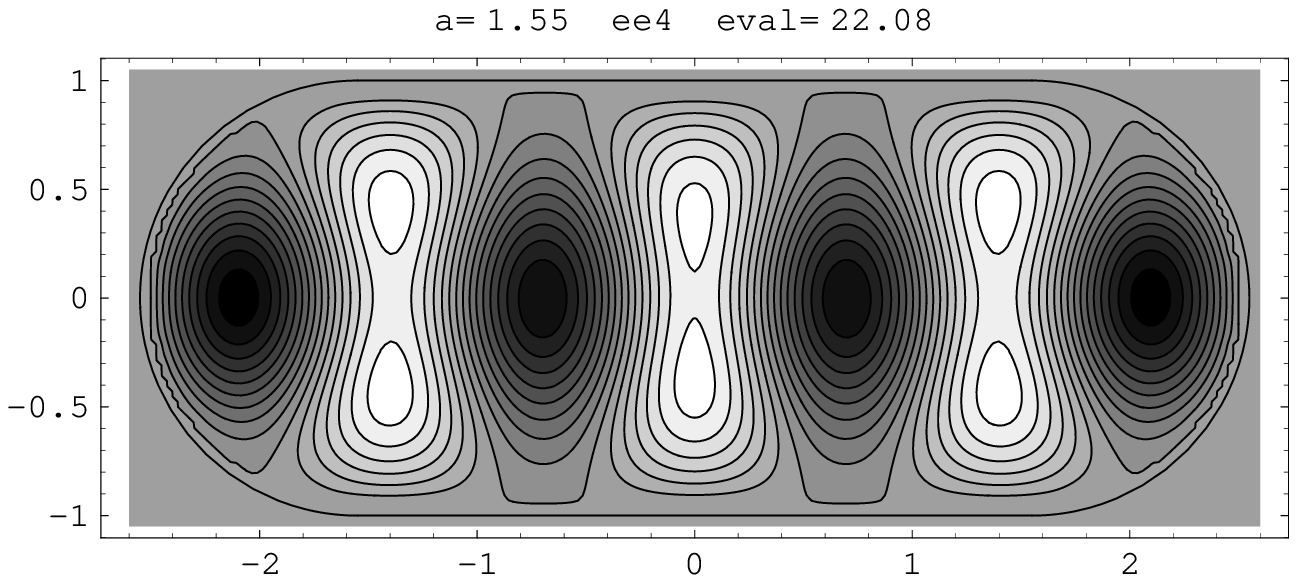}
\centereps{10.2cm}{4.cm}{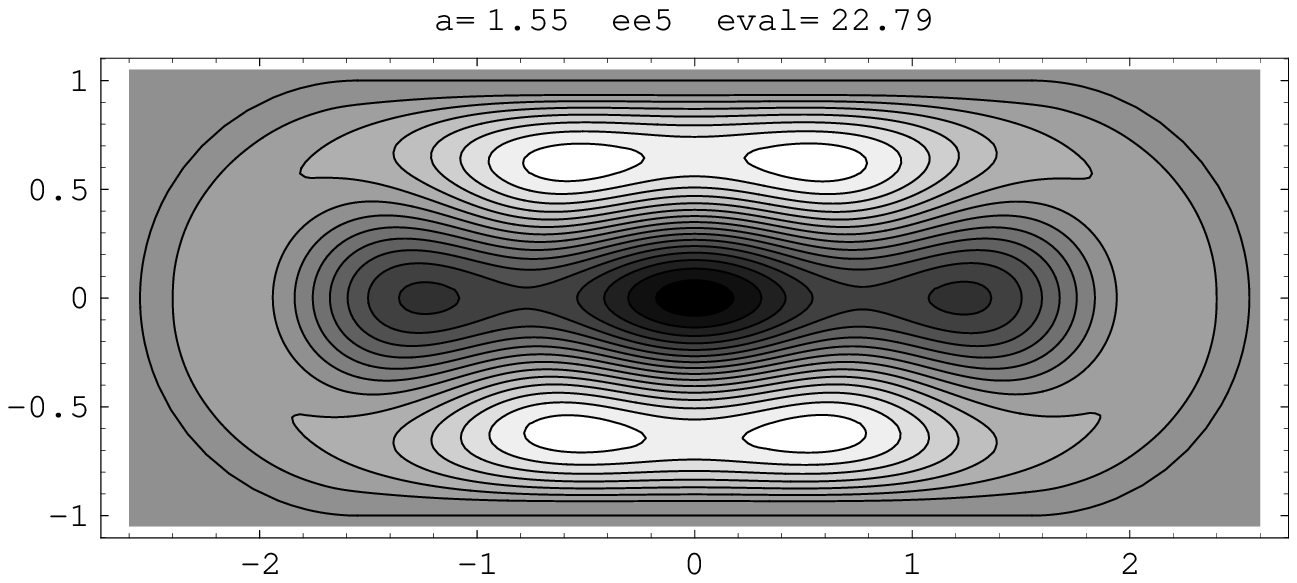}
}

\newpage{
\vspace*{.1in}
\centereps{10.8cm}{4.cm}{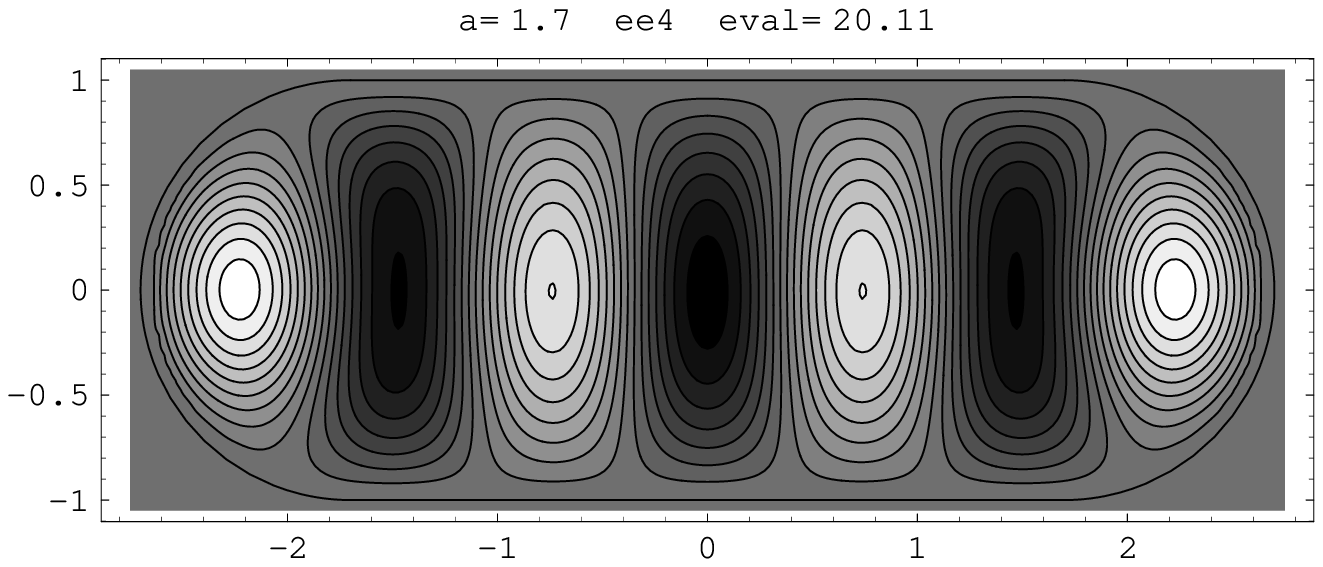}
\centereps{10.8cm}{4.cm}{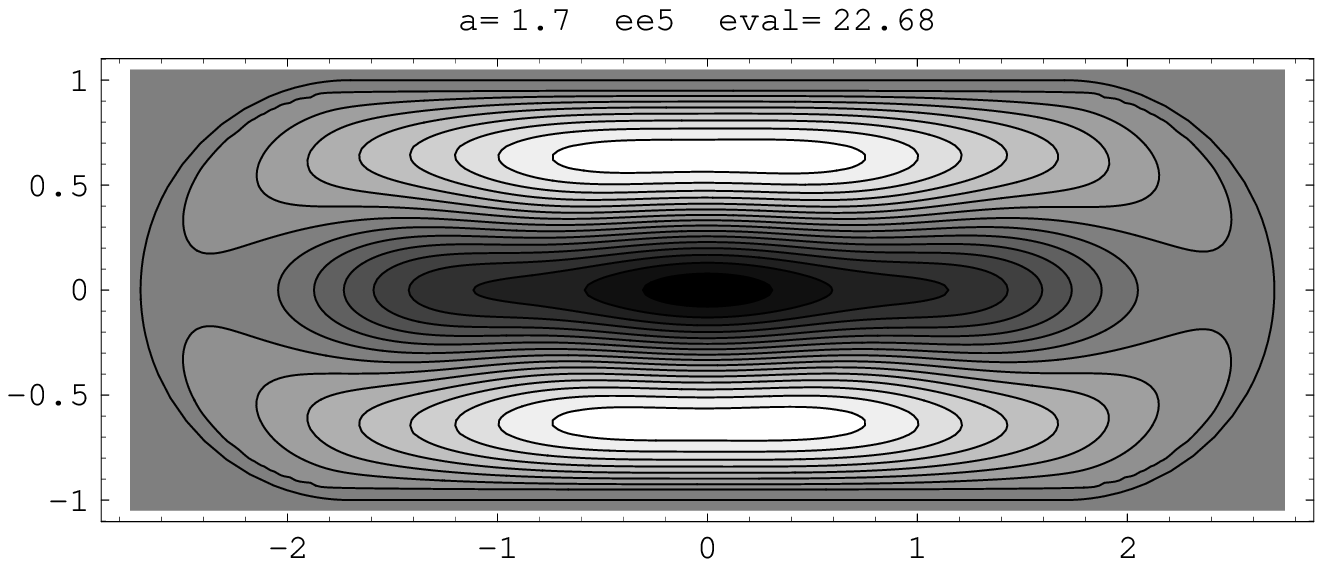}
\centereps{10.04cm}{4.cm}{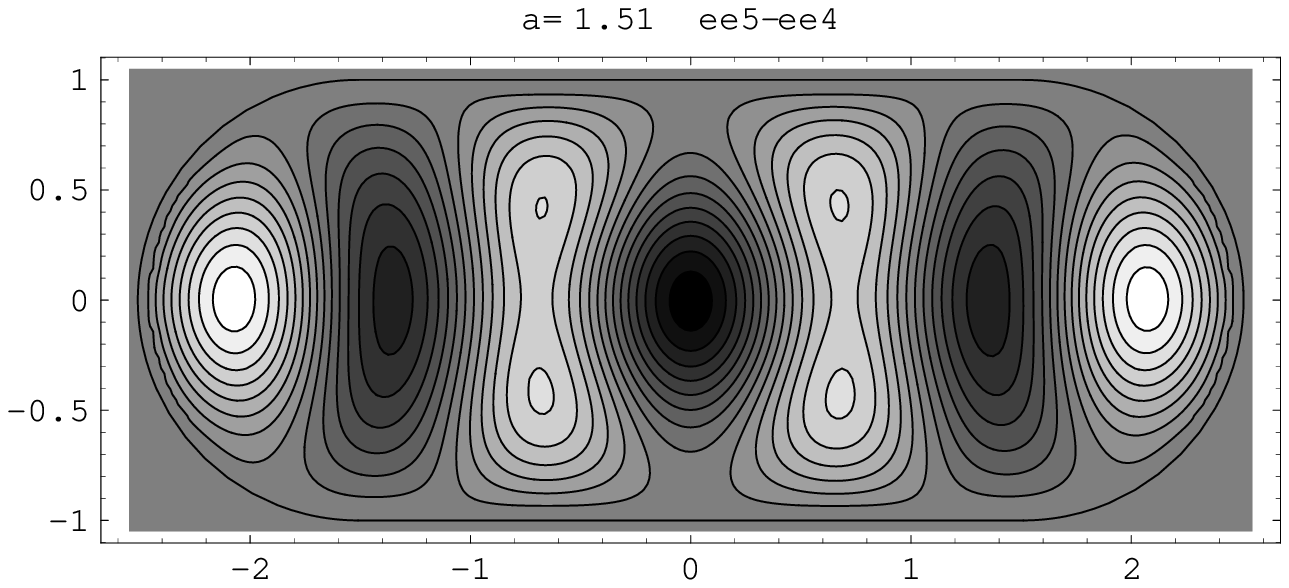}
\centereps{10.04cm}{4.cm}{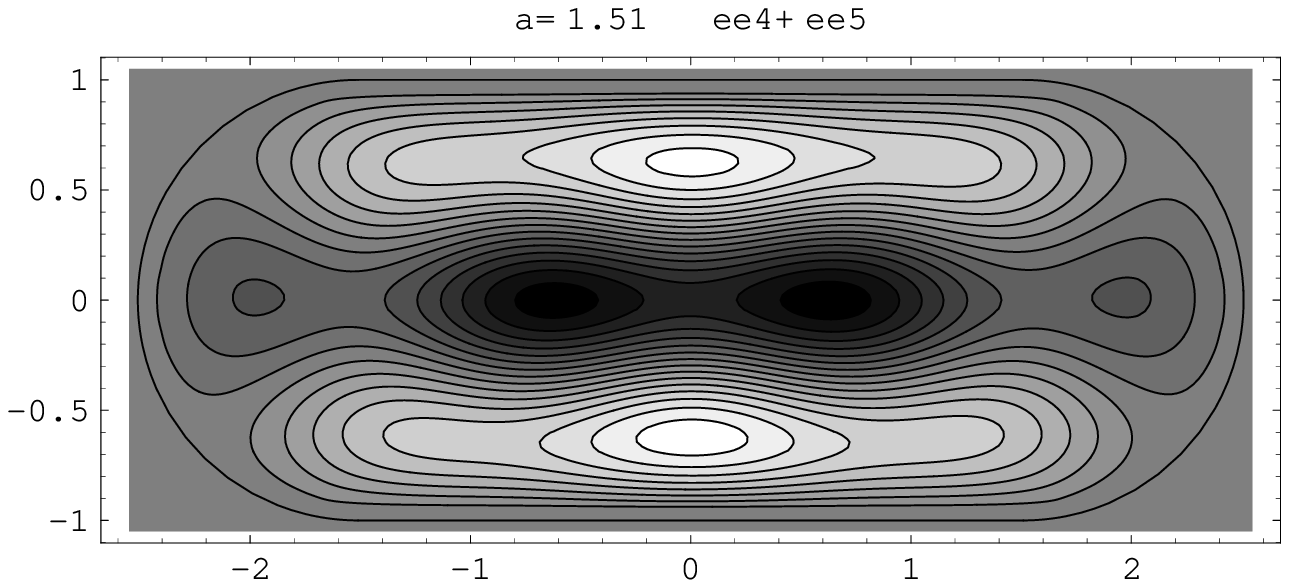}
}
 
\newpage

\noindent
\end{document}